\input amstex
\raggedbottom
\documentstyle{amsppt}
\magnification=1200
\pageheight{9.0 true in} \pagewidth{6.5 true in}
\pageno=1 \NoRunningHeads
\topmatter
\title
A note on the sums of powers of consecutive $q$-integers
\endtitle
\author
 YILMAZ  SIMSEK,   Dayeoul Kim,  Taekyun Kim,  and Seog-Hoon
Rim
\endauthor
\abstract
  In this paper we construct the $q$-analogue of Barnes's Bernoulli numbers
and polynomials of degree 2, which is an answer to a part of
Schlosser's question. Finally, we will treat the $q$-analogue of
the sums of powers of consecutive integers.
\endabstract
\thanks  2000 AMS Subject Classification:  11B68, 11S40.
\newline keywords and phrases :Sums of powers, Bernoulli Numbers, $q$-Bernoulli Numbers, zeta
function, Dirichlet series
\endthanks
\NoRunningHeads
\endtopmatter

\document

\head{ 1. Introduction }
\endhead

In 1713, J. Bernoulli first discovered the method which one can
produce those formulae for the sum $\sum_{j=1}^{n}j^{k}$, for any
natural numbers $k$ (cf. [1],[3],[6],[7],[15],[22]). The Bernoulli
numbers are among the most interesting and important number
sequences in mathematics. These numbers first appeared in the
posthumous work ``{ \it Ars Conjectandi''} (1713) by Jakob
Bernoulli(1654-1705) in connection with sums of powers of
consecutive integers ( Bernoulli(1713) or D. E. Smith(1959) see
[15]).

Let $q$ be an indeterminate which can be considered in complex
number field, and for any integer $k$ define the $q$-integer as
$$
[ k]_{q}=\frac{q^{k}-1}{q-1},\text{ (cf.
[11],[12],[13],[16],[17]).}
$$
Note that $\lim_{q\rightarrow 1}[k]_{q}=k.$ Recently, many authors
studied $q$-analogue of the sums of powers of consecutive
integers.

In [6], Garrett and Hummel gave a combinatorial proof of a
$q$-analogue of $\sum_{k=1}^{n}k^{3}=  {\binom {n+1}{2}}^{2}$ as
follows:
$$
\sum_{k=1}^{n}q^{k-1}\left( \frac{1-q^{k}}{1-q}\right) ^{2}\left(
\frac{ 1-q^{k-1}}{1-q^{2}}+\frac{1-q^{k+1}}{1-q^{2}}\right) ={
{n+1} \brack {2} }_{q}^{2}, $$ where
$$
 { {n} \brack {k} }_{q}=\prod_{j=1}^{k}\frac{1-q^{n+1-j}}{1-q^{j}}
$$
denotes the $q$-binomial coefficient.
 Garrett and Hummel, in their
paper, asked for a simpler $q$-analogue of the sum of cubes . As a
response to Garrett and Hummel's question, Warnaar gave a simple
$q$-analogue of the sum of cubes as follows:
$$
\eqalignno{& \sum_{k=1}^{n}q^{2n-2k}\frac{\left( 1-q^{k}\right)
^{2}\left( 1-q^{2k}\right) }{\left( 1-q\right) ^{2}\left(
1-q^{2}\right) }= { {n+1} \brack {2}}_{q}^{2}. &(1.1)}
$$

In [21], Schlosser took up on Garrett and Hummel's second
question.  Especially, he studied the $q$-analogues of the sums of
consecutive integers, squares, cubes, quarts and quints. He
obtained his results by employing specific identities for
very-well-poised basic hypergeometric series. However, Schlosser
could not find the $q$-analogue of the sums of powers of
consecutive integers of higher order, and left it as question. In
[15], T. Kim evaluated sums of powers of consecutive $q$-integers
as follows:

For any positive integers $n,k$($>1$), $h\in \Bbb{Z}$, let
$$
 S_{n,q^{h}}(k)=\sum_{j=0}^{n}q^{hj}[k]_{q}^{n}.$$

Then, he obtained the interesting formula for $S_{n,q}(k)$ below:

$$
S_{n,q}(k)=\frac{1}{n+1}\sum_{j=0}^{n}\binom{n+1}{j}
 \beta _{j,q}q^{kj}[k]_{q}^{n+1-j}-\frac{(1-q^{(n+1)k})\beta _{n+1,q}
}{n+1},
$$
where $\beta _{j,q}$ are the modified Carlitz's $q$-Bernoulli
numbers. Indeed, this formula is exactly a $q$-analogue of the
sums of powers of consecutive integers due to Bernoulli.

Recently, the problem of $q$-analogues of the sums of powers have
attracted the attention of several
authors([8],[9],[15],[19],[21],[23]). Let
 $$\eqalignno{&
S_{m,n}(q)=\sum_{k=1}^{n}[k]_{q^{2}}[k]_{q}^{m-1}q^{(n-k)\frac{m+1}{2}}.
&(1.2)}$$

 Then Warnaar[23] (for
$m=3$) and Schlosser[21] gave formulae for $m=1,2,3,4,5$ as the
meaning of the $q$-analogues of the sums of consecutive integers,
squares, cubes, quarts and quints. By two families of polynomials
and Vandermonde determinant, Guo and Zeng[9] found the formulae
for the $q$-analogues of the sums of consecutive integers (for
$m=1,2,...,5).$ They recovered the formulae of Warnaar and
Schlosser, for $\ m=6,7,...,11.$ In  [21],
 Schlosser speculated on the existence of a general formula for $
S_{m,n}(q),$ which is defined in (1.2), and left it as an open
problem.

By using T. Kim technical method to construct $q$-Bernoulli
numbers and polynomials in
[10],[11],[12],[13],[14],[15],[16],[19], we construct the
$q$-analogue of Barnes' Bernoulli numbers and polynomials of
degree 2, which is an answer to a part of Schlosser's question.
Finally, we give some formulae for $S_{m,n}(q)$.

We define generating function $F_{k,q}^{\ast }(t)$ of the
$q$-Bernoulli numbers $\beta _{n,k,q}^{\ast }$ ( $n\geq 0$ ) as
follows:
$$\eqalignno{
F_{k,q}^{\ast }(t) &=-t\sum_{j=0}^{\infty } q^{k-j}[j]_{q^{2}}\exp
(t[j]_{q} q^{\frac{k-j}{2}})  &(1.3) \cr &=\sum_{n=0}^{\infty }
\frac{\beta_{n,k,q}^{\ast }t^{n}}{n!} \text{ (cf. }
[10],[11],[12],[13],[14],[15],[16]).}$$

\proclaim{Theorem 1}
 Let $n,k$ be positive integers. Then
$$ \beta _{n,k,q}^{\ast }=\left(
\frac{1}{1-q}\right) ^{n-1}\sum_{m=0}^{n}\binom {n}{m}
 \frac{(-1)^{m}mq^{\frac{(n-1)(k-1)}{2}+k+m-2}}{(1-q^{m-\frac{n-1}{2}
-2})(1-q^{m-\frac{n-1}{2}})}.
$$
\endproclaim

We define generating function $F_{k,q}^{\ast }(t;k)$ of the
$q$-Bernoulli polynomials $\beta _{n,k,q}^{\ast }(k)$ ( $n\geq 0$
) as follows:
$$
F_{k,q}^{\ast }(t;k)=-t\sum_{j=0}^{\infty }q^{-j}[j+k]_{q^{2}}\exp
(t[j+k]_{q}q^{\frac{-j}{2}})=\sum_{n=0}^{\infty }\frac{\beta
_{n,k,q}^{\ast }(k)t^{n}}{n!}.$$

\proclaim{Theorem 2}
 Let $n,k$ be positive integers. Then
$$\beta _{n,k,q}^{\ast
}(k)=\frac{1}{[2]_{q}(1-q)^{n-2}}\sum_{m=0}^{n}\binom {n}{ m}
 (-1)^{m}\left( \frac{mq^{k(m-1)}}{1-q^{m-\frac{n-1}{2}-2}}-\frac{
mq^{k(m+1)}}{1-q^{m-\frac{n-1}{2}}}\right) .$$
\endproclaim

\proclaim{Theorem 3} ({\bf{General formula for }$S_{m,n}(q)$}) Let
$n,k$ be positive integers. Then
$$
S_{n,k}(q)=\sum_{j=0}^{k-1}[j]_{q^{2}}[j]_{q}^{n-1}q^{\frac{(n+1)(k-j)}{2}}=
\frac{\beta _{n,k,q}^{\ast }(k)-\beta _{n,k,q}^{\ast }}{n}.$$
\endproclaim

\head{ 2. Preliminary  }
\endhead

Let $n,k$ be positive integers ($k>1$), and let
$$
S_{n}(k)=\sum_{j=1}^{k-1}j^{n}=1^{n}+2^{n}+...+(k-1)^{n}.
$$
It is well-known that
$$\split
S_{1}(k) &=\frac{1}{2}k^{2}-\frac{1}{2}k, \cr S_{2}(k)
&=\frac{1}{3}k^{3}-\frac{1}{2}k^{2}+\frac{1}{6}k, \cr S_{3}(k)
&=\frac{1}{4}k^{4}-\frac{1}{2}k^{3}+\frac{1}{4}k^{2}, \cr
{}&\cdots .
\endsplit$$

Thus we have the following three conjectures:

(I) $S_{n}(k)$ is a polynomial in $k$ of degree $n+1$ with leading
coefficient $\frac{1}{n+1},$

(II) The constant term of $S_{n}(k)$\ is zero, i.e., $S_{n}(0)=0,$

(III) The coefficient of $k^{n}$ in $S_{n}(k)$ is $-\frac{1}{2}.$

Therefore, $S_{n}(k)$ is a polynomial in $k$ of the form

$$
S_{n}(k)=\frac{1}{n+1}k^{n+1}-\frac{1}{2}k^{n}+a_{n-1}k^{n-1}+\cdots+a_{1}k.
$$

We note that
$$
\frac{d}{dk}S_{n}(k)=k^{n}-\frac{n}{2}k^{n-1}+\cdots\text{ \ .}
$$

To make life easier, we put the first two conjectures together and
we reach the following conjecture, which is what Jacques Bernoulli
(1654-1705) claimed more than three hundred years ago.

{\bf{Bernoulli: }}There exists a unique monic polynomial of degree
$n$, say $B_{n}(x),$ such that
$$
S_{n}(k)=\sum_{j=1}^{k-1}j^{n}=1^{n}+2^{n}+...+(k-1)^{n}=
\int_{0}^{k}B_{n}(x)dx .
$$

As the $q$-analogue of $S_{n}(k)$, Schlosser[21]  considered the
existence of general formula on $S_{m,n}(q)$, and he gave the
below values:
$$
\sum_{k=1}^{n}[k]_{q^{2}}[k]_{q}^{m-1}q^{(n-k)\frac{m+1}{2}},
$$
where $m=1,2,...,5.$

Indeed, $$\eqalignno{
\sum_{k=1}^{n}[k]_{q^{2}}[k]_{q}q^{\frac{3}{2}(n-k)}&=\frac{
[n]_{q}[n+1]_{q}[n+\frac{1}{2}]_{q}}{[1]_{q}[2]_{q}[\frac{3}{2}]_{q}},
&\quad  m=2,\cr \sum_{k=1}^{n}[k]_{q^{2}}[k]_{q}^{2}q^{2(n-k)}&= {
{n+1} \brack {2} }^2_q ,&\quad m=3, \cr
\sum_{k=1}^{n}[k]_{q^{2}}[k]_{q}^{3}q^{\frac{5}{2}(n-k)} &=\frac{
(1-q^{n})(1-q^{n+1})(1-q^{n+\frac{1}{2}})}{(1-q)(1-q^{2})(1-q^{\frac{5}{2}})}
\cr &\times \left(
\frac{(1-q^{n})(1-q^{n+1})}{(1-q)^{2}}-\frac{q^{n}(1-q^{
\frac{1}{2}})}{1-q^{\frac{3}{2}}}\right) ,&\quad  m=4, \cr
\sum_{k=1}^{n}[k]_{q^{2}}[k]_{q}^{4}q^{3(n-k)} &=\frac{
(1-q^{n})^{2}(1-q^{n+1})^{2}}{(1-q)^{2}(1-q^{2})(1-q^{3})} \cr
&\times \left(
\frac{(1-q^{n})(1-q^{n+1})}{(1-q)^{2}}-\frac{q^{n}(1-q)}{
1-q^{2}}\right), &\quad m=5.}
$$

 T. Kim[15] proved the smart formula for the $q$-analogue of
$S_n(k)$ as follows:
$$
\sum_{k=0}^{n-1}q^{k}[k]_{q}=\frac{1}{2}\left(
[n]_{q}^{2}-\frac{[2n]_{q}}{ [2]_{q}}\right)$$ and
$$
\sum_{k=0}^{n-1}q^{k+1}[k]_{q}^{2}=\frac{1}{3}[n]_{q}^{3}-\frac{1}{2}\left(
[n]_{q}^{2}-\frac{[2n]_{q}}{[2]_{q}}\right)
-\frac{1}{3}\frac{[3n]_{q}}{ [3]_{q}}.
$$

\head{ 3. Proof of Main Theorems  }
\endhead

We define generating function $F_{k,q}^{\ast }(t)$ of the
$q$-Bernoulli numbers $\beta _{n,k,q}^{\ast }$ ( $n\geq 0$ ) as
follows:
$$\eqalignno{&
F_{k,q}^{\ast }(t)=-t\sum_{j=0}^{\infty }q^{k-j}[j]_{q^{2}}\exp
(t[j]_{q}q^{ \frac{k-j}{2}})=\sum_{n=0}^{\infty }\frac{\beta
_{n,k,q}^{\ast }t^{n}}{n!} .&(3.1)}$$

\demo{Proof of Theorem 1}
 Let
$$
\sum_{n=0}^{\infty }\frac{\beta _{n,k,q}^{\ast }t^{n}}{n!}
=-t\sum_{j=0}^{\infty }q^{k-j}[j]_{q^{2}}\exp
(t[j]_{q}q^{\frac{k-j}{2}}).
$$
By using Taylor series in the above then we have
$$
\sum_{n=0}^{\infty }\frac{\beta _{n,k,q}^{\ast }t^{n}}{n!}
=-t\sum_{j=0}^{\infty }q^{k-j}[j]_{q^{2}}\sum_{n=0}^{\infty
}\frac{ [j]_{q}^{n}q^{\frac{n(k-j)}{2}}}{n!}t^{n}.$$

By using some elementary calculations in the above, we have
$$\split
 \sum_{n=0}^{\infty }\frac{\beta _{n,k,q}^{\ast
}t^{n}}{n!} &= -t\sum_{j=0}^{\infty }q^{k-j}[j]_{q^{2}} \cr
&\times \sum_{n=0}^{\infty }\left\{ \left( \frac{1}{1-q}\right)
^{n}q^{ \frac{n(k-j)}{2}}\sum_{m=0}^{n}\binom{n}{m}
 (-1)^{m}q^{jm}\right\} \frac{t^{n}}{n!} \cr
&=\frac{-t}{q^{2}-1}\sum_{n=0}^{\infty }\left(
\frac{1}{1-q}\right)
^{n}q^{k+\frac{nk}{2}}\sum_{m=0}^{n}\binom{n}{m}
 (-1)^{m} \cr
&\times \sum_{j=0}^{\infty }\left(
q^{mj-j-\frac{jn}{2}}(1-q^{2j})\right) \frac{t^{n}}{n!} .
\endsplit$$

By using geometric power series in the above and after some
calculations, we obtain
$$\split
\sum_{n=0}^{\infty }\frac{\beta _{n,k,q}^{\ast }t^{n}}{n!}
&=\frac{-t}{ q^{2}-1}\sum_{n=0}^{\infty }\left(
\frac{1}{1-q}\right) ^{n} \cr &\times \sum_{m=0}^{n}
\binom{n}{m} \frac{(-1)^{m}q^{m+k+\frac{n}{2}(k-1)-1}(1-q^{2})}{\left( 1-q^{-1+m-%
\frac{n}{2}}\right) \left( 1-q^{1+m-\frac{n}{2}}\right)
}\frac{t^{n}}{n!}. \endsplit$$

For $n=0,$ then $\beta _{0,k,q}^{\ast }=0.$ Thus, we have
$$\split
\sum_{n=1}^{\infty }\frac{\beta _{n,k,q}^{\ast }t^{n}}{n!}
&=-t\sum_{n=1}^{\infty }\left( \frac{1}{1-q}\right) ^{n-1} \cr
&\times \sum_{m=1}^{n-1}\binom{n-1}{m-1}
 \frac{(-1)^{m-1}q^{m+k+\frac{(n-1)(k-1)}{2}-2}}{\left( 1-q^{-2+m-
\frac{n-1}{2}}\right) \left( 1-q^{m-\frac{n-1}{2}}\right)
}\frac{t^{n-1}}{ (n-1)!} \cr &=\sum_{n=1}^{\infty }\left(
\frac{1}{1-q}\right) ^{n-1} \cr
 &\times \sum_{m=0}^{n}\binom{n}{m} \frac{(-1)^{m}mq^{m+k+\frac{(n-1)(k-1)}{2}-2}}{\left( 1-q^{m-\frac{
n-1}{2}-2}\right) \left( 1-q^{m-\frac{n-1}{2}}\right)
}\frac{t^{n}}{n!}.
\endsplit$$
By comparing the coefficients of $\frac{t^{n}}{n!}$ on both sides
of the above equation, we easily arrive at the desired result.
\quad\quad\qed \enddemo

We define generating function $F_{k,q}^{\ast }(t;k)$ of the
$q$-Bernoulli polynomials $\beta _{n,k,q}^{\ast }(k)$ ( $n\geq 0$
) as follows:
$$\eqalignno{&
F_{k,q}^{\ast }(t;k)=-t\sum_{j=0}^{\infty }q^{-j}[j+k]_{q^{2}}\exp
(t[j+k]_{q}q^{\frac{-j}{2}})=\sum_{n=0}^{\infty }\frac{\beta
_{n,k,q}^{\ast }(k)t^{n}}{n!}. &(3.2)}$$

\demo{Proof of Theorem 2}
 Let
$$
\sum_{n=0}^{\infty }\frac{\beta _{n,k,q}^{\ast }(k)t^{n}}{n!}
=-t\sum_{j=0}^{\infty }q^{-j}[j+k]_{q^{2}}\exp
(t[j+k]_{q}q^{\frac{-j}{2}}).$$

By using Taylor expansion of $e^{x}$ in the above, we have
$$
\sum_{n=0}^{\infty }\frac{\beta _{n,k,q}^{\ast }(k)t^{n}}{n!}
=-t\sum_{j=0}^{\infty }q^{-j}[j+k]_{q^{2}}\sum_{n=0}^{\infty
}\frac{ [j+k]_{q}^{n}q^{\frac{-nj}{2}}}{n!}t^{n}.
$$

By using some elementary calculations in the above, we have
$$\split
\sum_{n=0}^{\infty }\frac{\beta _{n,k,q}^{\ast }(k)t^{n}}{n!}
&= -t\sum_{j=0}^{\infty }q^{-j}[j+k]_{q^{2}} \\
& \times \sum_{n=0}^{\infty }\left\{ \left( \frac{1}{1-q}\right)
^{n}\sum_{m=0}^{n}\binom{n}{m}
(-1)^{m}q^{m(j+k)-\frac{nj}{2}}\right\} \frac{t^{n}}{n!} \cr
=&\frac{t}{q^{2}-1}\sum_{n=0}^{\infty }\left( \frac{1}{1-q}\right)
^{n}\sum_{m=0}^{n}\binom{n}{m}
 (-1)^{m}q^{mk} \cr
&\times \sum_{j=0}^{\infty }\left( 1-q^{2(j+k)}\right)
q^{-j+jm-\frac{jn}{2} }\frac{t^{n}}{n!} .\endsplit$$

By using geometric power series in the above and after some
calculations, we obtain
$$\split
\sum_{n=0}^{\infty }\frac{\beta _{n,k,q}^{\ast }(k)t^{n}}{n!} &=
-t\sum_{n=0}^{\infty }\frac{1}{[2]_{q}}\left( \frac{1}{1-q}\right)
^{n-1} \cr &\times \sum_{m=0}^{n-1}\binom{n}{m}
 \frac{(-1)^{m+1}q^{\frac{n-1}{2}(k-1)+k+m-1}}{\left( 1-q^{-1+m-\frac{
n-1}{2}}\right) \left( 1-q^{1+m-\frac{n-1}{2}}\right)
}\frac{t^{n-1}}{(n-1)!}. \endsplit$$

Thus, we get
$$\split  \sum_{n=0}^{\infty }\frac{\beta
_{n,k,q}^{\ast }(k)t^{n}}{n!} &=-t\sum_{n=0}^{\infty
}\frac{1}{[2]_{q}}\left( \frac{1}{1-q}\right) ^{n-1} \cr
  &\times
\sum_{m=0}^{n}\binom{n}{m}
 (-1)^{m}\left( \frac{q^{mk}}{1-q^{m-1-\frac{n}{2}}}-\frac{q^{(m+2)k}
}{1-q^{m+1-\frac{n}{2}}}\right) \frac{t^{n}}{n!}\endsplit$$ and
$$\split
\sum_{n=1}^{\infty }\frac{\beta _{n,k,q}^{\ast }(k)t^{n}}{n!}
&=\sum_{n=1}^{\infty }\frac{1}{[2]_{q}}\left( \frac{1}{1-q}\right)
^{n-2} \cr &\times \sum_{m=0}^{n}\binom{n}{m} (-1)^{m}\left(
\frac{mq^{(m-1)k}}{1-q^{m-\frac{n-1}{2}-2}}-\frac{
mq^{(m+1)k}}{1-q^{m-\frac{n-1}{2}}}\right) \frac{t^{n}}{n!} .
\endsplit$$

By comparing the coefficients of $\frac{t^{n}}{n!}$ on both sides
of the above equations, we easily arrive at the desired result.
\quad\quad\qed \enddemo

\demo{Proof of Theorem 3} Let
$$\eqalignno{
&-\sum_{j=0}^{\infty }q^{-j}[j+k]_{q^{2}}\exp
(t[j+k]_{q}q^{\frac{-j}{2} })+\sum_{j=0}^{\infty
}q^{k-j}[j]_{q^{2}}\exp (t[j]_{q}q^{\frac{k-j}{2}}) &(3.3) \cr
&=\sum_{j=0}^{k-1}q^{k-j}[j]_{q^{2}}\exp
(t[j]_{q}q^{\frac{k-j}{2}}).}$$

By using Theorem 1 and Theorem 2 and (3.3), we easily arrive at
the desired result. Our proof of Theorem 3 runs parallel that of
Theorem 1 and Theorem 2 above, so we choose to omit the details
involved.
\enddemo

\remark{Remark 1}  The Barnes double zeta function is defined by
$$
\zeta _{2}(s,w\mid w_{1},w_{2})=\sum_{m,n=0}^{\infty
}(w+mw_{1}+nw_{2})^{-s}, \ \ \text{{Re}}(s)>2, \text{ (cf.
[2])},$$ for complex number $w\neq 0,w_{1},w_{2}$ with positive
real parts.

The Barnes' polynomial: \ Barnes[2] introduced $r$-tuple Bernoulli
polynomials $_{r}S_{m}(u;\widetilde{w})$ by,
$$\eqalignno{ &
F_{r}(t;u;\widetilde{w})=\frac{(-1)^{r}te^{-ut}}{
\prod_{j=1}^{r}(1-e^{-w_{j}t})}=\sum_{k=1}^{r}\frac{_{r}S_{1}^{(k+1)}(u;
\widetilde{w})(-1)^{k}}{t^{k-1}}&(3.4)\cr &+\sum_{m=1}^{\infty
}\frac{ _{r}S_{m}^{^{\prime
}}(u;\widetilde{w})}{m!}(-1)^{m-1}t^{m}, }
$$
for $\mid t\mid <\min \left\{ \frac{2\pi }{\mid w_{1}\mid
},...,\frac{2\pi }{ \mid w_{r}\mid }\right\} $.  Here
$w_{1},w_{2},...,w_{r}$ are complex number with positive real
parts,  $\widetilde{w}=(w_{1},w_{2},...,w_{r})$ and
$_{r}S_{1}^{^{(k)}}(u;\widetilde{w})$  the $k$-th derivative of
$_{r}S_{1}(u; \widetilde{w})$ with respect to $u.$

Substituting $u=-1,\ w_{1}=w_{2}=-1$ and $r=2$ into (3.4), we have
$$
F_{2}(t;-1;-1,-1)=\frac{te^{t}}{(1-e^{t})^{2}}.
$$
In (3.1), $$
 \lim_{q\rightarrow 1}F_{k,q}^{\ast }(t)
=\lim_{q\rightarrow 1}-t\sum_{j=0}^{\infty }q^{k-j}[j]_{q^{2}}\exp
(t[j]_{q}q^{\frac{k-j}{2} })=-F_{2}(t;-1;-1,-1).$$
\endremark

\vskip 20pt
 \head{4. Further Remarks and Observations on a class
of $q$- Zeta Functions}\endhead

In this section, by using generating functions of $F_{k,q}^{\ast
}(t)$ and $ F_{k,q}^{\ast }(t;k)$, we produce new definitions of
$q$-polynomials and numbers. These generating functions are very
important in  the case of defining $ q $-zeta function. Therefore,
by using these generating functions and Mellin transformation, we
will define  the $q$-zeta function.

By applying Mellin transformation in (3.2), we obtain
$$
\frac{1}{\Gamma (s)}\int_{0}^{\infty }t^{s-2}F_{k,q}^{\ast
}(-t;k)dt=\sum_{n=0}^{\infty
}\frac{[n+k]_{q^{2}}q^{-n\frac{s+2}{2}}}{ [n+k]_{q}^{s}},
$$
where $\Gamma (s)$  denotes the Euler gamma function.

 For $s\in \Bbb C$, we define
$$
\zeta _{k,q}^{\ast }(s;k)=\sum_{n=0}^{\infty
}\frac{[n+k]_{q^{2}}q^{-n\frac{ s+2}{2}}}{[n+k]_{q}^{s}},\quad
\text{Re}(s)>2.$$

By Mellin transformation in (3.1), we obtain
$$\split
\frac{1}{\Gamma (s)}\int_{0}^{\infty }t^{s-2}F_{k,q}^{\ast }(-t)dt
&=\sum_{n=0}^{\infty
}\frac{[n]_{q^{2}}q^{\frac{(k-n)(2-s)}{2}}}{[n]_{q}^{s} } \cr
&=\zeta _{k,q}^{\ast }(s),\text{Re}(s)>2.\endsplit$$

For any positive integer $n,$ Cauchy Residue Theorem in the above
equation, we have
$$
\zeta _{k,q}^{\ast }(1-n)=-\frac{\beta _{n,k,q}^{\ast }}{n}.
$$

\vskip 20pt \quad  YILMAZ  SIMSEK

\quad Mersin University, Faculty of Science,

\quad Department of Mathematics 33343 Mersin, Turkey

\quad e-mail: ysimsek$\@$mersin.edu.tr

\vskip 20pt

\quad Daeyeoul Kim

\quad Department of Mathematics,

\quad Chonbuk National University,

\quad Chonju, 561-756, Korea

\quad   e-mail: dykim$\@$math.chonbuk.ac.kr

\vskip 20pt \quad  TAEKYUN KIM

\quad Institute of Science Education,

\quad  Kongju National University Kongju 314-701, Korea

\quad  e-mail: tkim$\@$kongju.ac.kr

\vskip 20pt

\quad  Seog-Hoon Rim

\quad Department of Mathematics Education,

\quad  Kyungpook National University, Taegu 702-701,  Korea

\quad  e-mail : shrim$\@$knu.ac.kr

    \enddocument